\documentclass[12pt,a4paper,reqno]{amsart}
 
\usepackage[margin=1in]{geometry}
\usepackage{amsmath,amsthm,amssymb}
\usepackage{colortbl}
\usepackage{xcolor}
\usepackage{graphicx}
\usepackage{subcaption}
\usepackage{lscape}
\usepackage{rotating}
\usepackage{multirow}
\usepackage{adjustbox}

\usepackage{setspace}
\usepackage{graphicx}

\usepackage{hyperref}
\usepackage[square,numbers]{natbib}
\usepackage{notoccite}
\usepackage{comment}
\usepackage{enumitem, kantlipsum}

\begin{document}

\title{An Improved Autoencoder Conjugacy Network to Learn Chaotic Maps}

\author[]{Meagan Carney}
\address{The University of Queensland School of Mathematics and Physics} \email{m.carney@uq.edu.au}
\urladdr{}

\author[]{Cecilia Gonz\'{a}lez-Tokman}
\address{The University of Queensland School of Mathematics and Physics} \email{cecilia.gt@uq.edu.au}
\urladdr{}

\author[]{Ruethaichanok Kardkasem}
\address{The University of Queensland School of Mathematics and Physics} \email{r.kardkasem@uq.edu.au}
\urladdr{}

\author[]{Hongkun Zhang}
\address{Great Bay University, China} \email{zhanghk@gbu.edu}
\urladdr{}

\thanks{MC and HZ would like to thank the Erwin Schr\"{o}dinger Institute, Vienna, Austria, where part of this work was completed. CGT, MC and HZ would like to thank the University of Queensland Raybould Fellowship for support on this project. MC and RK would like to thank Great Bay University, China, where part of this work was completed.}

\begin{abstract}
We introduce a method for learning chaotic maps using an improved autoencoder neural network that incorporates a conjugacy layer in the latent space. The added conjugacy layer transforms nonlinear maps into a simple piecewise linear map (the tent map) whilst enforcing dynamical principles of well-known and defective conjugacy functions that increase the accuracy and stability of the learned solution. We demonstrate the method's effectiveness on both continuous and piecewise chaotic one-dimensional maps and numerically illustrate improved performance over related traditional and recently emerged deep learning architectures.
\end{abstract}

\maketitle

\section{Introduction}
One of the greatest challenges in working with chaotic dynamical systems from a computational perspective is generating reliable orbits to make predictions and inform theoretical proofs under the constraints of an exponentially increasing error. A common solution is to use some advanced numerical method to decrease the error in the system as much as possible and inform these results using theoretical insights on small perturbations of chaotic systems \cite[Chapter 9]{evtbook}. With the recent advancements in deep learning, it is natural to ask whether these methods can be used to improve our understanding of chaotic dynamical systems. Many studies have already started examining the possibility of this approach \cite{chattopadhyay2020data, boulle2020classification, fan2020long, haluszczynski2021controlling, jin2016modeling, lusch2018deep, pathak2018model, serrano2021new, wang2019neural, weng2022modeling}.

With the introduction of Poincar\'{e} sections to analyse complex behaviour of dynamical systems by transforming their continuous evolution into discrete 'snapshots', chaotic maps have become fundamental in capturing both chaotic and deterministic dynamics within complex systems. These maps are particularly useful in revealing the sensitive dependence on initial conditions and the bifurcations that arise under shifting parameters. Highly accurate numerical modelling of chaotic maps can provide insights into the long-term and statistical behaviour of a dynamical system, when closed-form solutions cannot be readily derived.

\begin{figure}[htbp]
    \centering
    \begin{subfigure}[b]{\textwidth}
            \includegraphics[width=\textwidth]{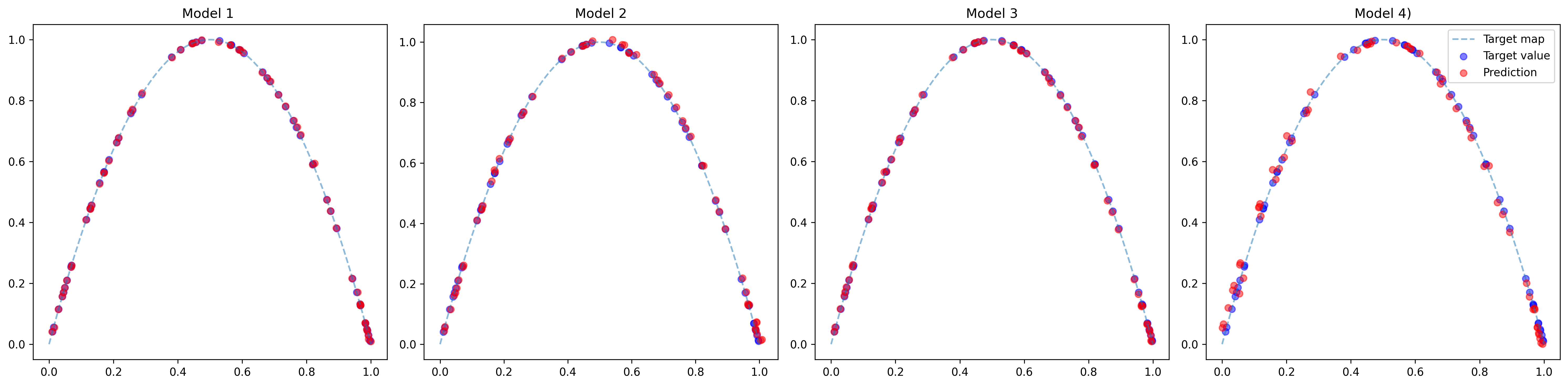}
            \caption{Step size = 2}
            \label{fig:logis_step_2}
    \end{subfigure}
    \begin{subfigure}[b]{\textwidth}
            \includegraphics[width=\textwidth]{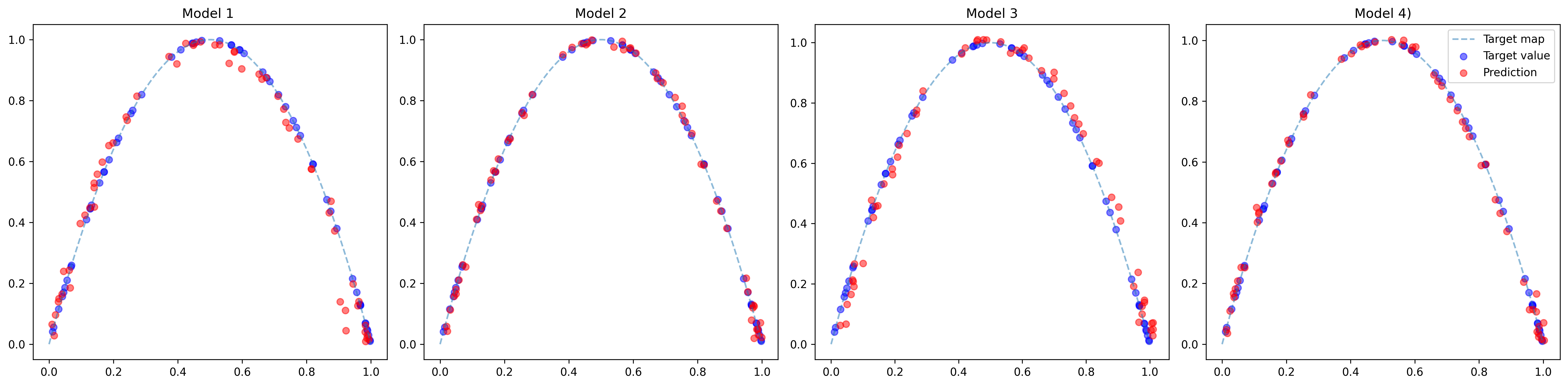}
            \caption{Step size = 5}
            \label{fig:logis_step_5}
    \end{subfigure}
        \begin{subfigure}[b]{\textwidth}
            \includegraphics[width=\textwidth]{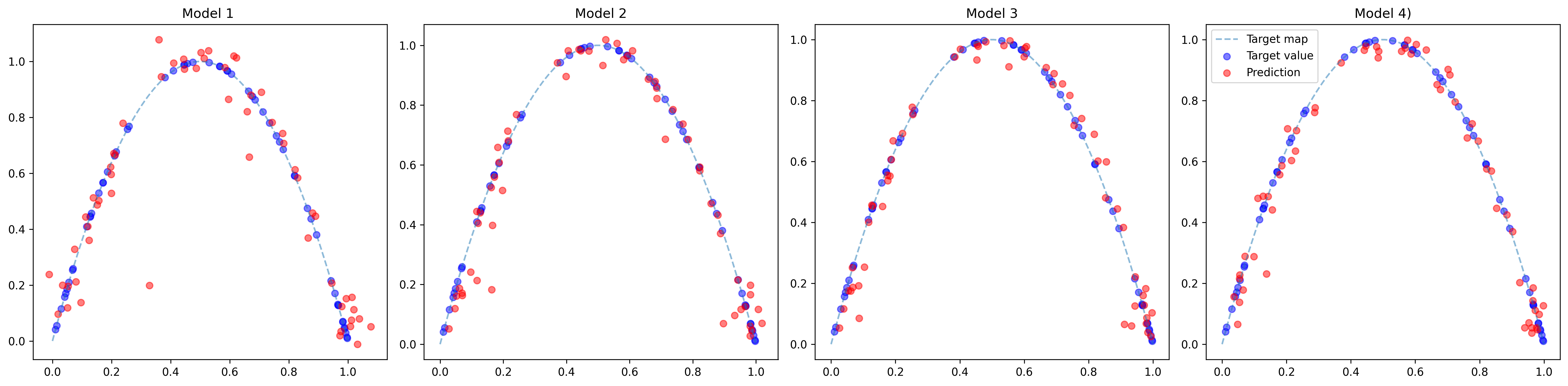}
            \caption{Step size = 7}
            \label{fig:logis_step_7}
    \end{subfigure}
    \caption{The prediction of the test dataset using deep learning models in \ref{auto}, including AENN (Models 1 and 2), FNN (Model 3), and PINN (Model 4), is based on logistic map ($r=4, x_0=0.4$) with time-series inputs of varying step window sizes. The models used sequences of 300 in length, with 80\% of the data allocated for training and 20\% for testing.}
    \label{fig:logis_step}
\end{figure}

From this lens, the goal of a neural network is to learn a single step of a chaotic map to the highest possible precision; however, the way by which this is learned varies depending on the structure of the network. For example, in the data driven feed-forward neural network (FNN), the true numerical error is dependent on our ability to trust generated orbital steps (which we cannot if the system has sensitive dependence on initial conditions), to represent the true values along the orbit of a target chaotic map upon which the network is trained. One possible approach is to borrow ideas from Physics Informed Neural Networks (PINNs), originally designed to solve partial differential equations (PDEs) \cite{raissi2019physics}. A PINN is a supervised neural network that uses the functional representation of a PDE and a set of data generated from the PDE to learn and forward predict values. The benefit of a PINN is that it is robust for lower quantities of data which makes it a possible candidate to learn a chaotic map where long, accurate orbits cannot be generated. However, recursively, we are faced with the difficulty of generating an orbit accurate enough, even for small stretches, to serve as a training set for the PINN. In fact, here we find that traditional methods of numerically generating long orbits to estimate statistical properties of even simple maps, such as piecewise expanding maps —often involving the addition of low order stepwise noise (on the order of \(10^{-6}\)) and the application of the shadowing lemma — requires careful (and often computationally intensive) calibration of the hyper-parameters and often results in a failure to converge to the true map (see figure \ref{fig:logis_step} for an illustrative example). Although not a main focus of this investigation, such a result emphasises the need for mathematicians to contribute to the literature surrounding the stability and reliability of deep learning for chaotic systems.

Given the issues surrounding orbit generation for training a deep network, we turn our attention to the class of unsupervised deep networks known as \textit{autoencoders}. A network is called an autoencoder if its input is trained to copy its output. Consequently, the benefit of using an autoencoder is that one can train the network to learn the chaotic map by using the first time step of a large uniform sample of initial conditions, completely circumventing the need to use a long orbit from a single initial condition.

Autoencoder network architectures have become widely used for learning dynamical systems, albeit often from a more applied, rather than a theoretical perspective. Before introducing some of the applied work that has inspired this investigation, we would like to emphasise that there exists a rich field of work that has recently emerged in the space of dynamical systems and deep networks. We refer the reader to \cite{PhysicaDspecialissue} for a nice archive of related work in this area, including, but not specific to autoencoder approaches. 

In general, an autoencoder will take in feature inputs (in this case, the points of our chaotic map), compress the inputs into the hidden layer called the \textit{latent space} where it performs coefficient estimation on a nonlinear combination of activation functions to learn the function describing the map and finally decompresses the data back into its original form as an output. Consequently, the interpretability of the resulting learned model is difficult (or impossible) when autoencoders are employed in their general form.

In the applied dynamical systems space, the authors in \cite{otto2019linearly} presented a combination of an autoencoder and a linear recurrent network with numerical case studies on a handful of nonlinear flows while \cite{agostini2020exploration, eivazi2020deep, lee2020model,  xu2020multi, simpson2021machine} combined autoencoders with other machine learning algorithms, such as long short-term memory (LSTM) networks and convolutional neural networks (CNN), to develop reduced-order models for nonlinear flow modelling. The authors in \cite{yang2022learning} use the autoencoder structure with physical principles included in the latent space (similar to the PINN) for modelling well-studied differential equations in engineering. The authors in \cite{baig2023autoencoder} showed that an autoencoder can be used to compress high-dimensional data while eliminating extraneous noise without losing information. These results were in the context of drug-resistant pathogen emergence informed by systems of growth dynamics. 

A common theme emerges in the applied autoencoder literature: whilst deep learning models are often regarded as black boxes, one can build a deep network architecture based on known dynamical properties that can aid in interpreting the learned representation (or output) of the underlying dynamical system. From a dynamics viewpoint, coordinate transformations arise as natural techniques for building dynamically informed autoencoder architectures. This is because the objective of a coordinate transformation in dynamics is to study the behaviour of a highly complex system by reducing it to a simpler system where conclusions can be drawn. 

Common coordinate transformations utilised in dynamics include the Koopman operator, for transforming a nonlinear system into a linear one in an infinite-dimensional space by acting on observable functions, allowing for linear analysis techniques; Poincar\'{e} sections to reduce a continuous system to a lower-dimensional discrete map by recording intersections with a chosen surface, making it easier to reveal long term behaviour of the analogous complex system; and the conjugacy for establishing equivalencies in dynamics from complex properties to systems with more easily understood dynamical or geometric properties. These ideas have been interwoven with deep network architecture in recent literature. For example, \cite{lusch2018deep} introduced the Koopman operator in the latent space to transform certain nonlinear dynamical systems into their linear counterparts. This approach allows for learning the linear mapping and subsequently reversing the process to generate the orbits of the nonlinear system. Similarly, \cite{champion2019data} developed a parsimonious (SINDy) representation for coordinate transformation within an autoencoder network. Recent work by \cite{bramburger2021deep} leverages Poincar\'{e} sections and conjugate mappings (or generally, commuter mappings \cite{skufca2007relaxing, skufca2008concept}) to create an autoencoder architecture that outputs higher precision orbits of higher-dimensional chaotic systems by imposing a conjugate mapping to a lower dimensional chaotic map in the latent space. Whilst this model works well for systems that would be computationally expensive to run on their own, in certain settings it can suffer from numerical issues surrounding the approximation of an unknown conjugacy mapping in the latent space. In particular, there is a requirement of the latent space (conjugacy function) to take the nonlinear form $L(y) = c_1y+c_2y^2$ where $c_1$ and $c_2$ are learned parameters, sensitive to initial choices and suffer from gradient vanishing (see Table \ref{tab:loss} for an experimental comparison). We therefore aim to improve on this autoencoder architecture by introducing a layer in the latent space that converts nonlinear maps into a simple piecewise linear map via an additional conjugate mapping, imposing a layer of stability and allowing for more accurate approximations. We illustrate these claims on numerical case studies of one-dimensional chaotic maps.

The outline of this manuscript is as follows:
\begin{itemize}
\item Section \ref{conj} describes the concept of a topological conjugacy and its use in the deep network architecture.
\item Section \ref{auto} introduces our autoencoder architecture which exploits a known conjugacy between the logistic map and the tent map to increase stability of the deep network.
\item Section \ref{exp} runs a series of numerical experiments using our architecture to learn a variety of one-dimensional chaotic maps and compares these results to known relevant architectures from literature.
\item Section \ref{conc} discusses the results and limitations of our investigation in context.
\end{itemize}

\section{Dynamical conjugacies, map geometry, and their roles in deep networks}\label{conj}

The notion of topological conjugacy provides a way to identify dynamical systems that differ only by a change of coordinates. Topological conjugacies are often exploited as a means of understanding a complex system by replacing it with a simplified one that has "nicer" properties, where proofs can be more easily managed. Under this framework, two maps $f: X\rightarrow X$ and $g: Y\rightarrow Y$ are termed \textit{topologically conjugate} if there exists some homeomorphism $h:X\rightarrow Y$, called a conjugacy function, such that

$$f\circ h = h\circ g.$$

\begin{figure}[ht]
    \centering
    \includegraphics[width=0.25\linewidth]{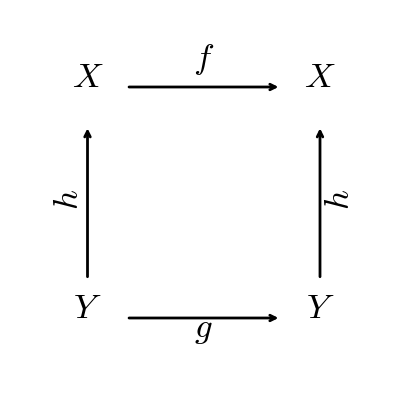}
    \caption{Commutative diagram}
    \label{fig:commutative_diagram}
\end{figure}

In the mathematical framework, the benefit of proving a conjugacy for some complex chaotic system is that it allows us to prove certain dynamical properties, such as ergodicity, invariance, and limiting distributions, for the simplified system and then apply the conjugacy to immediately recover these results in the more complex setting. In a similar way, deep networks aim to learn and approximate complex mappings by transforming them into more manageable representations, leveraging structured optimisation to extract and learn the essential features of a chaotic map. It is therefore the goal of this investigation to exploit known equivalencies in dynamics (e.g. the conjugacy) to instruct the network on how these more manageable representations should be built.

The general act of "learning" in a deep network comes from optimising the weights on a combination (possibly nonlinear) of basis functions (called activation functions) to approximate, as closely as possible the data (or map) the network is attempting to learn. The width of the network (or number of layers) indicates the number of convolutions used to generate the output of the network, while the depth of a layer indicates how many weighted basis functions, of a particular type, are being used to approximate a single function inside the convolution. Some dynamical maps are geometrically simpler and hence easier for a network to learn. For example, the tent map (\ref{eq:tent}), can be much more easily approximated by a linear combination of common basis functions, such as ReLU or SeLU, than its conjugate, continuous counterpart - the logistic map (\ref{eq:logistic}). Hence, the nonlinear nature of the logistic map can increase rounding errors, requiring more careful handling during numerical simulations. 

It is natural then to consider network architectures that exploit the properties of conjugacy functions in their convolution process. That is, using the representation $f = \hat{h}^{-1} \circ g \circ \hat{h}(x)$ to learn $f$ from known function $g$ by approximating $h$ (by $\hat{h}$) and $h^{-1}$ (by $\hat{h}^{-1}$) on separate layers of the network (e.g. the \textit{latent space}). Indeed, the authors in \cite{bramburger2021deep} illustrate that this idea yields very successful results for certain chaotic systems. However, there remains some ambiguity regarding the resulting trajectory generated using this method. For example, the function $h$ used to generate the orbit of the chaotic map $f$ is an approximation of the true conjugacy function and hence, is itself a defective conjugacy, a term coined by the authors of \cite{bollt2010comparing}, generating a map with the approximate dynamics of $f$ where some orbits between $f$ and $g$ may no longer be shared via the defective conjugacy. Adding known dynamical criterion to the latent space of the network defined by \cite{bramburger2021deep} helps impose a more rigorous dynamical equivalence and hence, a more stable output of the chaotic map $f$. In the sections to follow, we illustrate that this approximation of $h$ and $h^{-1}$ can result in unreliable orbits for certain one-dimensional expanding maps and propose the addition of a stability layer that relies on the analytic (true) conjugacy function between the tent and the logistic map.

\subsection{Tent map and Logistic map}
\subsubsection*{Tent Map}
The tent map is defined by the piecewise function:
    \begin{equation} \label{eq:tent}
        T(x_n) = x_{n+1} =
            \begin{cases}
            \mu \, x_n,  &x_n < 0.5 \\
            \mu (1 - x_n),  &x_n \geq 0.5
            \end{cases}
    \end{equation}
where \( \mu \) is a parameter.

\subsubsection*{Logistic Map} 
The logistic map is defined by the function:
    \begin{equation}\label{eq:logistic}
        L(x_n) = x_{n+1} = r \, x_n (1 - x_n)
    \end{equation}
where \( r \) is a parameter.

\subsubsection*{Conjugacy Between Tent Map and Logistic Map}
A conjugacy between the chaotic logistic map with parameter $r=4$ and the chaotic tent map with parameter $\mu=2$ is given by:
    \begin{equation} \label{eq:conjugacy}
        \phi(x) = \frac{2}{\pi} \arcsin\left(\sqrt{x}\right)
    \end{equation}
where \( x \in [0, 1] \). The inverse of the conjugacy function \( \phi(x) \) is given by:
    \begin{equation} \label{eq:conjugacy_inverse}
        \phi^{-1}(x) = \sin^2\left(\frac{\pi x}{2}\right).
    \end{equation}
Hence, the conjugacy function $\phi(x)$ is continuous, bijective and has a continuous inverse on the interval $[0,1]$. Thus, it is a homeomorphism.

This use of tent and logistic maps arise in various applications, such as chaos-based cryptography, image encryption, and pseudo-random sequence generators. \cite{xu2012theorem, guo2020quadratic, xu2020multi} Consequently, understanding equivalencies that can be exploited for ease and accuracy of deep network training of related chaotic maps, such as the maps illustrated here, is of particular importance to computer science and security.

\section{Autoencoder Neural Network}\label{auto}

\subsection{Network Architecture}
We introduce an Autoencoder Neural Network (AE) designed to embed a learned conjugacy mapping between a one-dimensional chaotic map of interest and the logistic map within the network’s latent space, following the approach of \cite{bramburger2021deep}. However, a key distinction in our architecture is the explicit enforcement of the known conjugacy mapping between the logistic map and the tent map in the latent space, thereby incorporating additional structural constraints into the network.

    \begin{figure}
        \centering
        \includegraphics[width=1\linewidth]{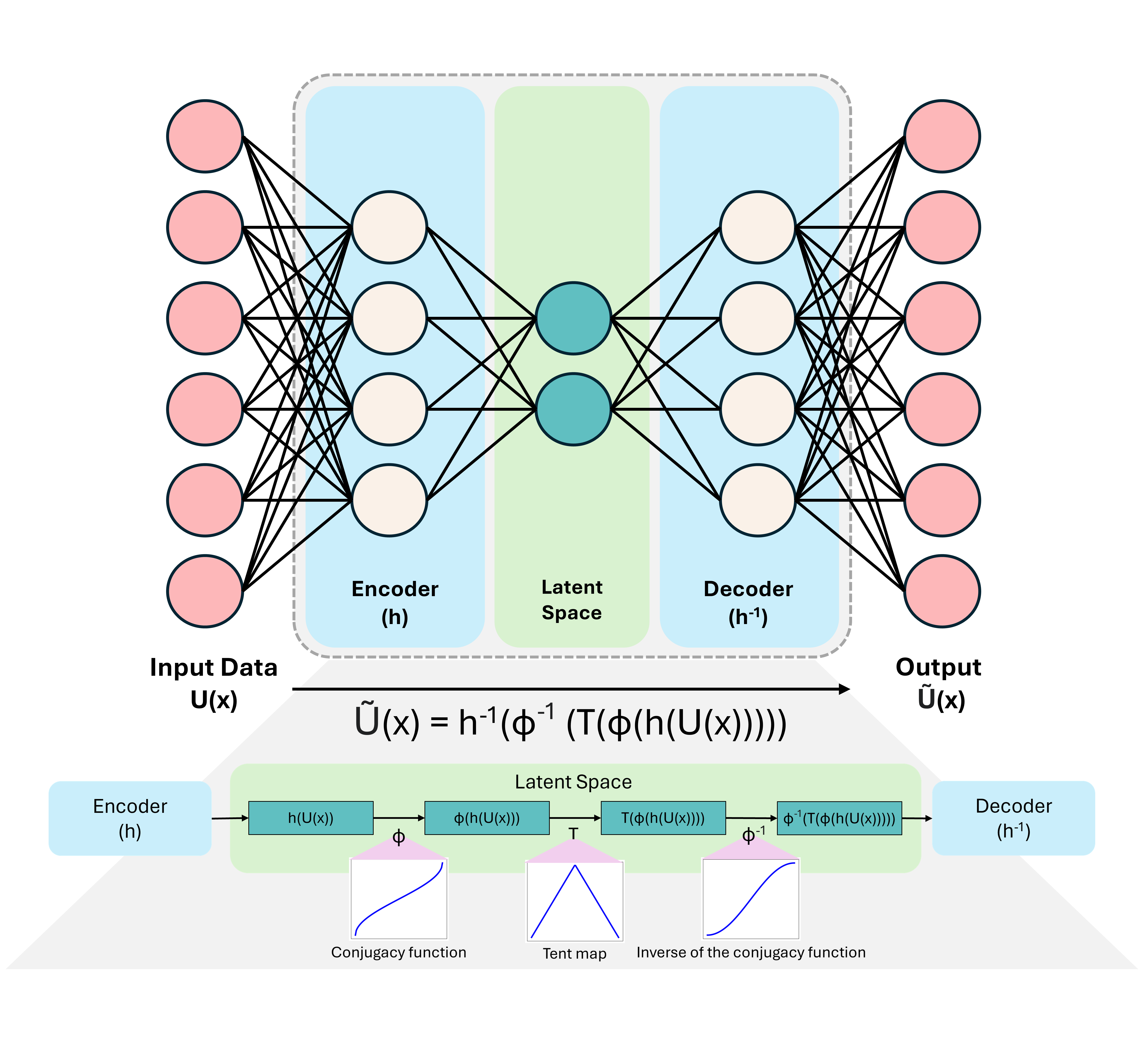}
        \caption{Diagram of the Autoencoder Neural Network architecture for finding a nonlinear map $U(x)$. The network first encodes the input into a latent space by applying the learned conjugacy function $h$. Then, it performs a coordinate transformation using the known conjugacy $(\phi)$ of the chaotic tent map $(T)$ and logistic map $(L)$. Finally, the unknown map $U(x)$ is recovered by applying the learned inverse conjugacy function $h^{-1}$.}
        \label{fig:AENN}
    \end{figure}

    \begin{figure}
        \centering
        \includegraphics[width=0.5\linewidth]{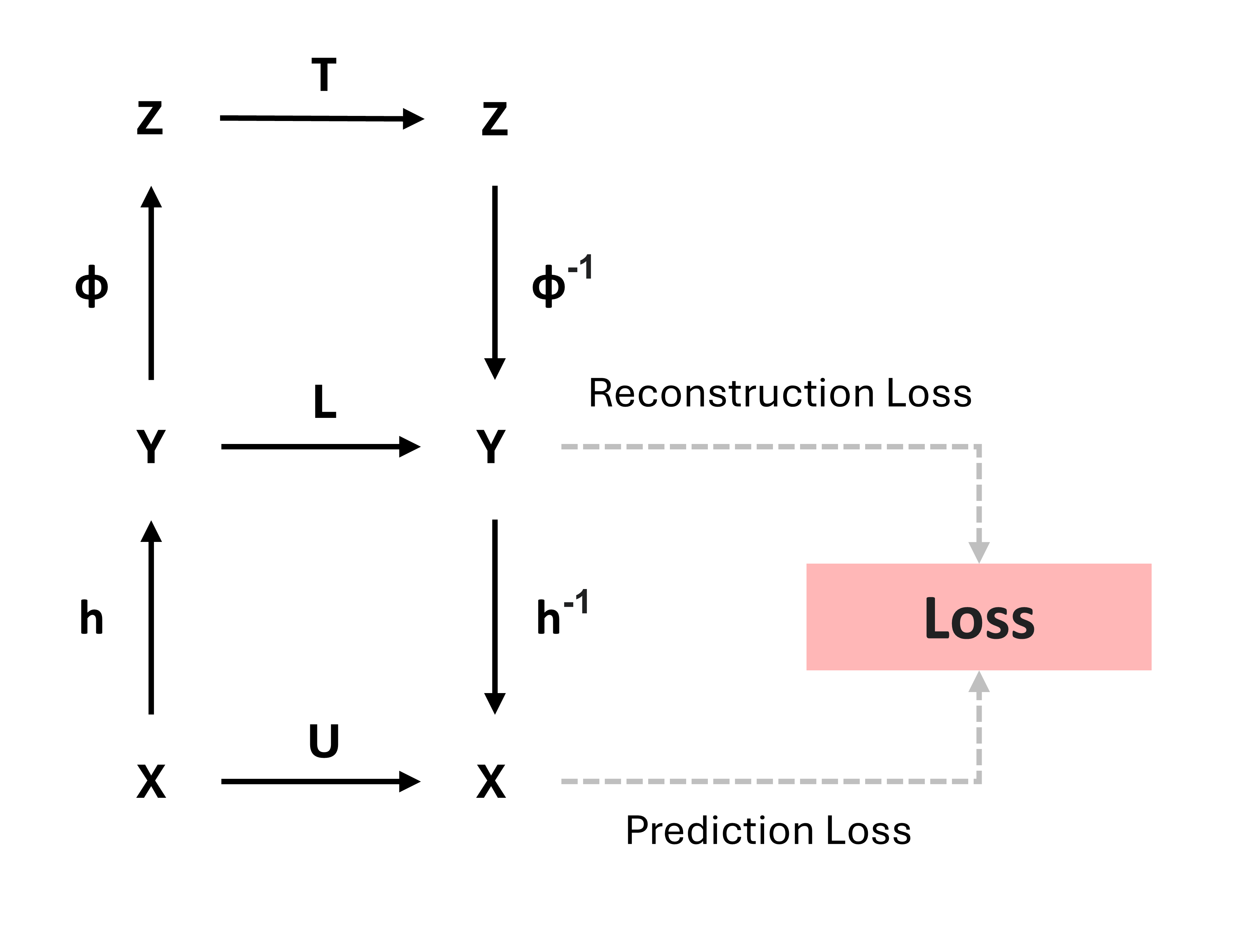}
        \caption{Commutative diagram of the AE embedding conjugacy function and tent map function in the latent space. Instead of attempting to predict an unknown nonlinear map $U(x)$ through the logistic map $(L)$, we use the conjugacy $(\phi)$ and the conjugate mapping of the logistic map, which is the tent map $(T)$. Therefore, we have two losses in the proposed AE: the reconstruction loss (\ref{eq:loss_recon}) and the prediction loss (\ref{eq:loss_pred})} 
        \label{fig:aenn_1}
    \end{figure}

\subsection{Latent space}
We begin with a subset of uniformly sampled points $\mathbf{x}$ from our one-dimensional map of interest. In the latent space, the data $\mathbf{x}$ is passed through an encoder $h$ to obtain the nonlinear variable $\mathbf{y}$. Next, we define the known transformation $\phi$ to take $\mathbf{y}$ from the latent space to an intermediate space within which the tent map is applied, and finally apply the known transformation inverse $\phi^{-1}$ to compute a new latent variable $y_{n+1} = \phi^{-1}(T(\phi(y_n))$, where $T$ is the tent map. The new latent variable is passed through a decoder $h^{-1}$, also in the latent space of the AE network, to obtain the final approximation of our map of interest. In a perfect approximation of the network, the variable $\mathbf{y}$ is expected to represent points along the logistic map.

Activation functions are used to approximate the encoder and decoder functions of the network. The encoder transforms the input data $\mathbf{x}$ into a latent representation $\mathbf{y}$ using the function \( h \) and an activation function \( \sigma \):
    \begin{equation}
        \mathbf{y} = h(\mathbf{x}) = \sigma(\mathbf{W}\mathbf{x} + \mathbf{b}),
    \end{equation}
when \( \mathbf{W} \) is weight matrix of the encoder, \( \mathbf{b} \) is  bias vector of the encoder. The decoder reconstructs the encoded data from the latent space using the function \( h^{-1} \), structured similarly to the encoder:
    \begin{equation}
        \Tilde{\mathbf{x}} = h^{-1}(\mathbf{y}) = \sigma(\mathbf{W}'\mathbf{y} + \mathbf{b}'),
    \end{equation}
when \( \mathbf{W}' \) is the weight matrix of the decoder (not necessarily the exact inverse of \( \mathbf{W} \), but it conceptually serves for input reconstruction), and \( \mathbf{b}' \) is the bias vector of the decoder.

Our model is an autoencoder neural network based on a feedforward architecture that predicts the data in the following forms: 
    \begin{align}
        \Tilde{U}(x) &  = h^{-1} \circ L \circ h, \\
        L(y) & = \phi^{-1} \circ T \circ \phi, \\
        \Tilde{U}(x) & = h^{-1} \circ \phi^{-1} \circ T \circ \phi \circ h
    \end{align}
where $T$ is the tent map in (\ref{eq:tent}), $L$ is the logistic map, (\ref{eq:logistic}) $\phi$ is the known conjugacy function in (\ref{eq:conjugacy}), $\phi^{-1}$ is inverse of the known conjugacy function in (\ref{eq:conjugacy_inverse}), $h$ is unknown (network approximated) encoder and $h^{-1}$ is unknown (network approximated) decoder. 
    
\subsection{Loss function }

The network incorporates both reconstruction and prediction mechanisms. The model is structured around an encoder-decoder architecture, coupled with a transformation and the conjugacy between the chaotic tent map and the logistic map in the latent space. The total loss function is composed of two main terms: a reconstruction loss that encourages identity mapping, and a prediction loss that captures the dynamics of a target function, as shown in Figure \ref{fig:aenn_1}. Let $U(x)$ be the target map, then the loss function is the sum total of the following loss:

    \begin{itemize}
        \item \textbf{Reconstruction Loss}: The reconstruction loss penalises deviation between the input \( \mathbf{x} \) and its reconstruction from the latent space through the decoder:
            \begin{equation} \label{eq:loss_recon}
                \mathcal{L}_{\text{recon}} =  \left\| U(x) - h^{-1}(h(U(x))) \right\|_{\text{MSE}}
            \end{equation}
        \item \textbf{Prediction Loss}: The prediction loss evaluates how well the transformed latent vector, when decoded, approximates the true target \( \mathbf{y} \):
            \begin{equation} \label{eq:loss_pred}
                \mathcal{L}_{\text{pred}} = \left\| U(x) - h^{-1}(\phi^{-1}(T(\phi(h(U(x)))))) \right\|_{\text{MSE}} 
            \end{equation}
    \end{itemize}

\section{Numerical Experiments}\label{exp}
We generate data from various one-dimensional chaotic maps to assess the predictive capabilities of our proposed model. In particular, we report results for the following piecewise and continuous one-dimensional chaotic maps:

\subsection{Custom Map}
The custom map is defined on the interval $[0,1)$ by the equation:
    \begin{equation}
        x_{n+1} = 16 \, x_n \left(1 - 2 \sqrt{x_n} + x_n\right).
    \end{equation}
The custom map features a right-skewed bell curve, which has been shown to be conjugate to the logistic map according to \cite{ruslan2022behavior}.

\subsection{Katsura-Fukuda Map}
The Katsura-Fukuda map is defined on the interval $[0,1)$ by the equation:
    \begin{equation}
        x_{n+1} = \frac{4 \, x_n (1 - x_n) (1 - r \, x_n)}{(1 - r \, x_n^2)^2}
    \end{equation}
where \( r \) is a parameter within the interval \( (0, 1) \). The map is a continuous map in the form of a regularised logistic map. This map becomes equivalent to the logistic map when $r=0.5$ \cite{glendinning2022differentiable}, recognised for its intricate bifurcation patterns and chaotic dynamics.

\subsection{Doubling Map}
The doubling map is defined on the interval $[0,1)$ by the equation:
    \begin{equation}
        x_{n+1} = 2x_n \mod 1, 
    \end{equation}
The doubling map, also known as the Bernoulli shift map, is a simple yet strongly chaotic map. The doubling map and the logistic map are topologically semi-conjugate \cite{layek2015introduction}.

\subsection{Pomeau-Manneville Map}
The Pomeau-Manneville map is defined on the interval $[0,1)$ by the equation:
    \begin{equation}
       x_{n+1} = x_n + ax_n^z \mod 1,
    \end{equation}
where \( z  > 1 \) is a parameter that controls the nonlinear behaviour. The map is commonly used to study intermittency in dynamical systems, combining regular and chaotic behaviour, making it useful for analysing transitions between different dynamical regimes. The Bernoulli map is a special case of this, with the parameter $a=z=1$ \cite{nee2018survival}.

\begin{figure}
    \centering
    \includegraphics[width=0.75\linewidth]{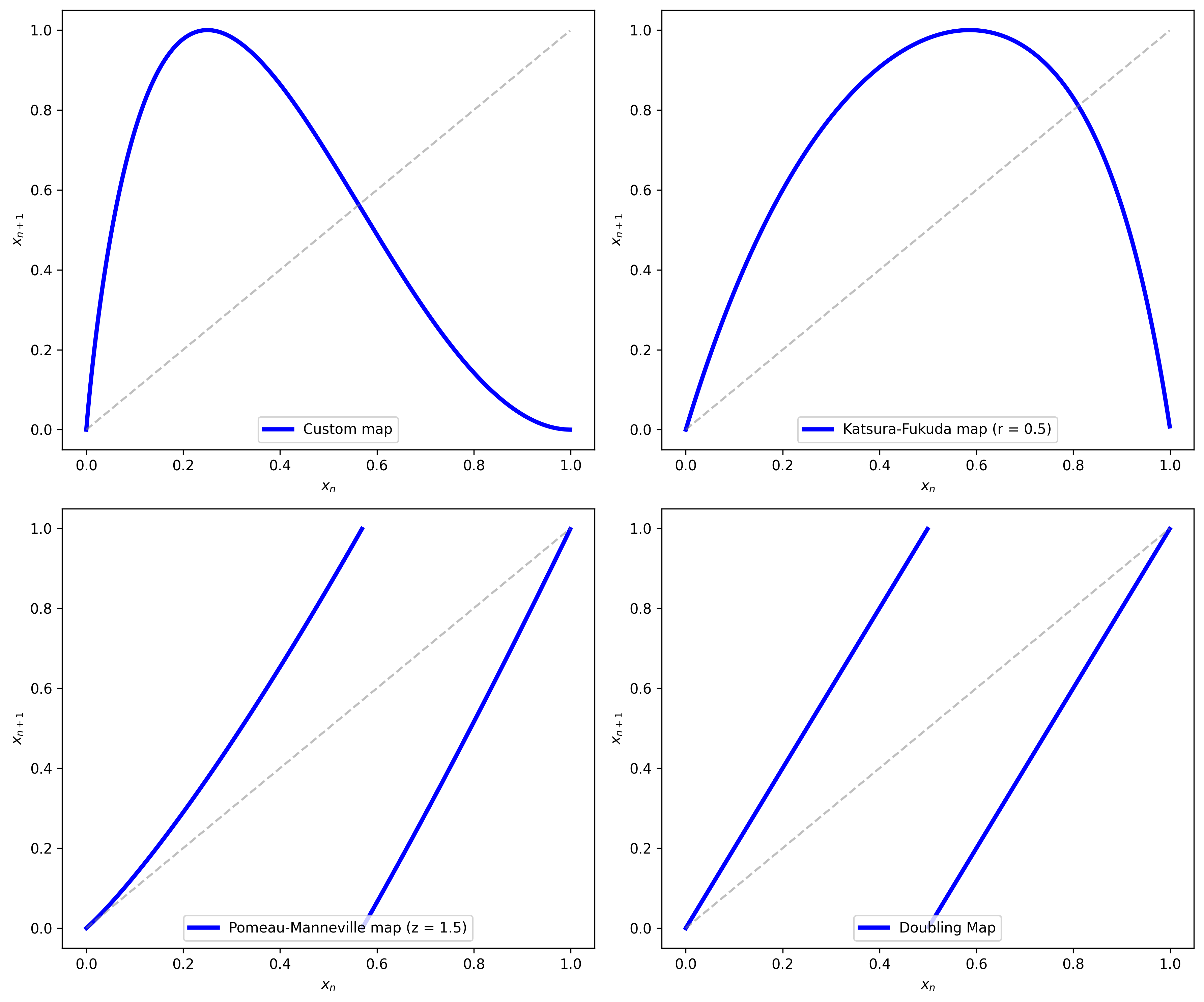}
    \caption{Custom map (top left), Katsura-Fukuda map (top right), Pomeau-Manneville map (bottom left) and Doubling map (bottom right)}
    \label{fig:target_map}
\end{figure}

For comparison, we analyse our proposed network (Model 1) against the comparable autoencoder structure in \cite{bramburger2021deep} (Model 2), the standard neural network (Model 3), and the physics-informed neural network (Model 4). We described the structure of these models for our equivalent purpose in detail below.
    \begin{itemize}
        \item \textbf{Model 1} Autoencoder Neural Network embedding a conjugacy between tent map and logistic map in the latent space. Our model includes the transformation of encoded vector by using the conjugacy between the logistic map (\ref{eq:logistic}) with parameter $r=4$ and tent map (\ref{eq:tent}) with parameter $\mu=2$.
        \item \textbf{Model 2} Autoencoder Neural Network embedding the logistic map in the latent space. This model is based on the network architecture from \cite{bramburger2021deep}, which the latent space includes the logistic map $L(y) = c_1y + c_2y^2$ with learnable parameters $c_1$ and $c_2$.
        \item \textbf{Model 3} Feedforward Neural Network (FNN): This network predicts the dynamical system by learning a purely data-driven mapping from past states to future states. This model employs supervised deep learning, where the network is trained on a dataset generated from the system.
        \item \textbf{Model 4} Physics-informed Neural Network (PINN): This network directly integrates the governing equation of the chaotic map into the training process by minimising data mismatches, while the residuals of the equations ensure that the network aligns with known dynamics \cite{raissi2019physics}.
    \end{itemize}
    
Various hyperparameters were tuned to train the networks and optimise model performance effectively. For continuous maps, we used the SeLU activation function, while the ReLU function is employed for piecewise maps. We implemented two hidden layers for continuous maps—one layer for the encoder and one for the decoder in the AE models. However, increasing the hidden layers is required for piecewise maps to enhance model performances. Details of these hyperparameter settings are provided in the Appendix.

\subsection{Confidence Interval}

Uncertainty in deep learning predictions is important for making reliable results and enhancing their interpretation, which is particularly useful for decision-making in certain applications. In this study, the confidence interval represents the prediction uncertainties. We use Monte Carlo (MC) dropout \cite{srivastava2014dropout, gal2016dropout} and ensemble methods \cite{lakshminarayanan2017simple} to estimate the confidence interval of predictions. Both methods provide multiple predictions for a given input, and the variation among these predictions reflects model uncertainty. The MC dropout technique estimates uncertainty by performing multiple stochastic forward passes with randomly dropped neurons at a specified probability within a single model, allowing the network to run multiple times with different dropout configurations each time. This method aims to reduce over-fitting by preventing specific neurons from becoming overly reliant on particular features. The mean is computed by averaging multiple forward passes through the same model with dropout applied, expressed as:
    \begin{equation} \label{eq:mcdropout}
        \mu_{\text{MC}} = \frac{1}{T} \sum_{i=1}^{N} \hat{x}^{(i)},
    \end{equation}
where \( T \) is the number of forward passes with dropout applied, and \( \hat{x}^{(i)} \) is the prediction from the \( i \)-th forward pass. 

Ensemble methods use multiple models on a single dataset through different initialisations or architectures. Each model in the ensemble is trained independently, and their predictions are combined to estimate the mean expressed as:
    \begin{equation} \label{eq:essem}
        \mu_{\text{Ensemble}} = \frac{1}{M} \sum_{j=1}^{M} \hat{x}^{(j)},
    \end{equation}
where \( M \) is the number of models in the ensemble, and \( \hat{x}^{(j)} \) is the prediction from the \( j \)-th model.

\section{Results and Discussion}\label{conc}
\subsection{Results and Conclusions}

\begin{figure}[htbp]
    \centering
    \begin{subfigure}[b]{\textwidth}
        \includegraphics[width=\textwidth]{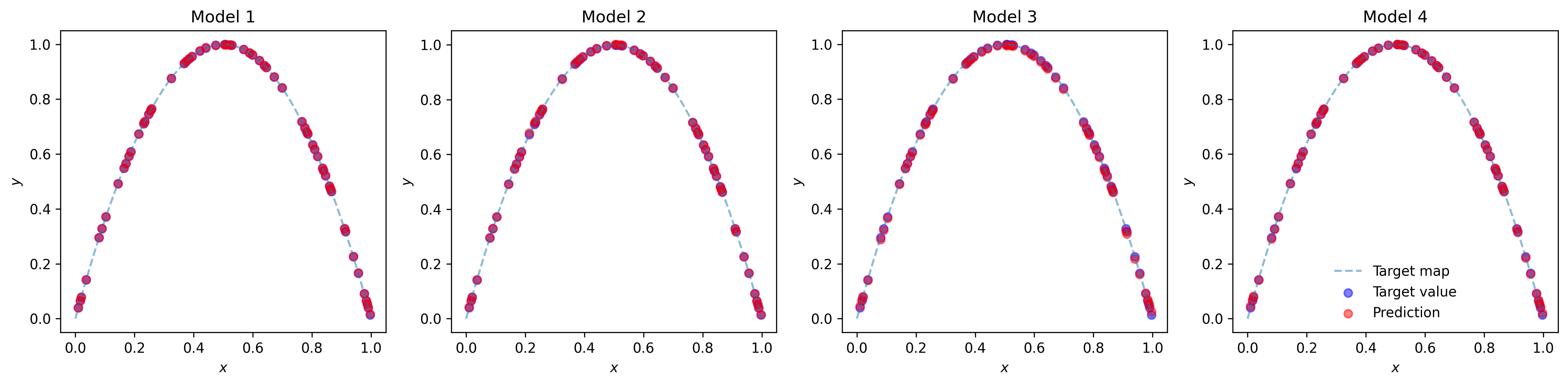}
        \caption{Logistic map }
        \label{fig:logistic_pred}
    \end{subfigure}

    \begin{subfigure}[b]{\textwidth}
        \includegraphics[width=\textwidth]{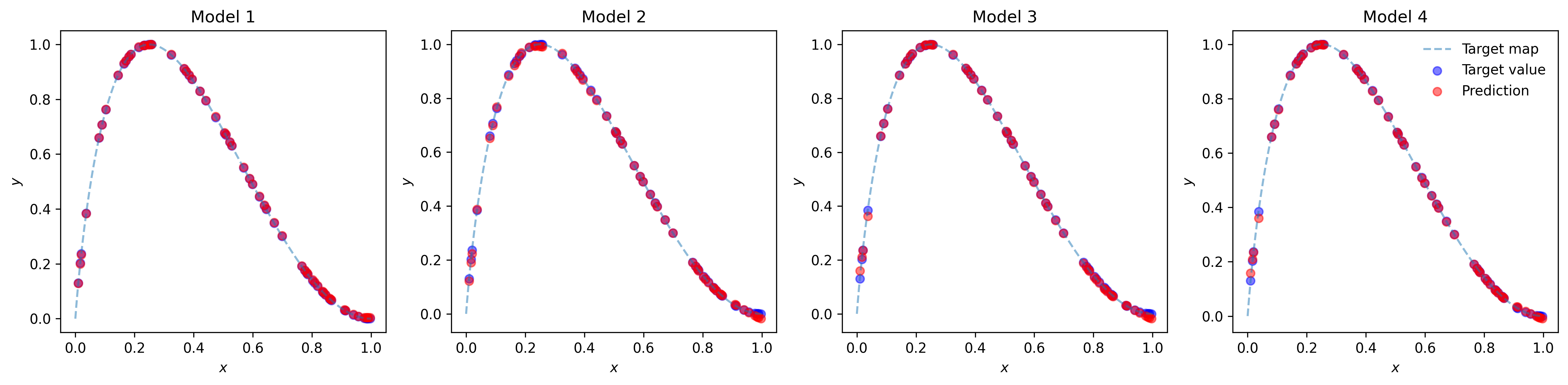}
        \caption{Custom map}
        \label{fig:custom_pred}
    \end{subfigure}
        
    \begin{subfigure}[b]{\textwidth}
        \includegraphics[width=\textwidth]{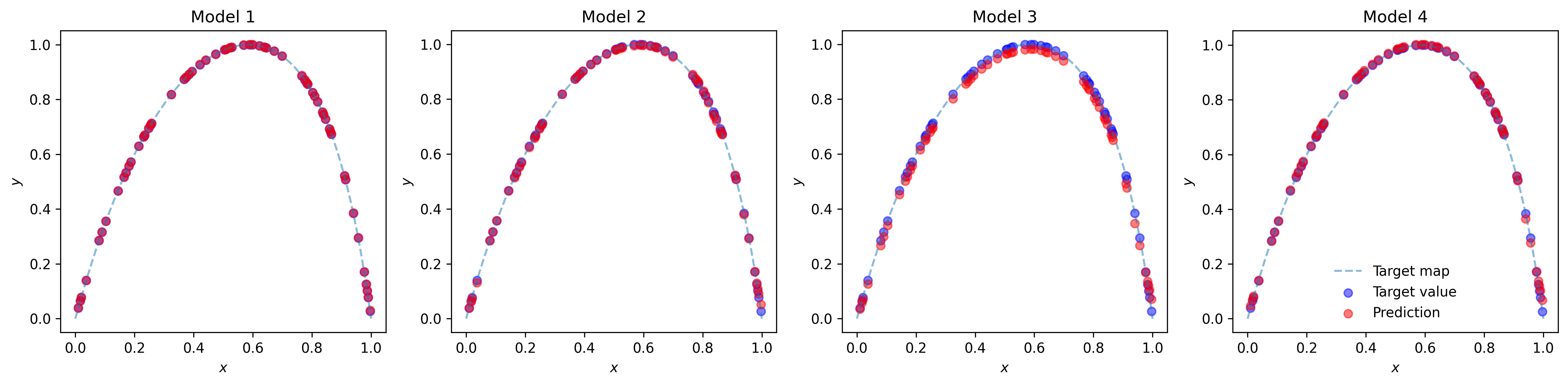}
        \caption{Katsura-Fukuda map}
        \label{fig:ks_pred}
    \end{subfigure}

    \begin{subfigure}[b]{\textwidth}
        \includegraphics[width=\textwidth]{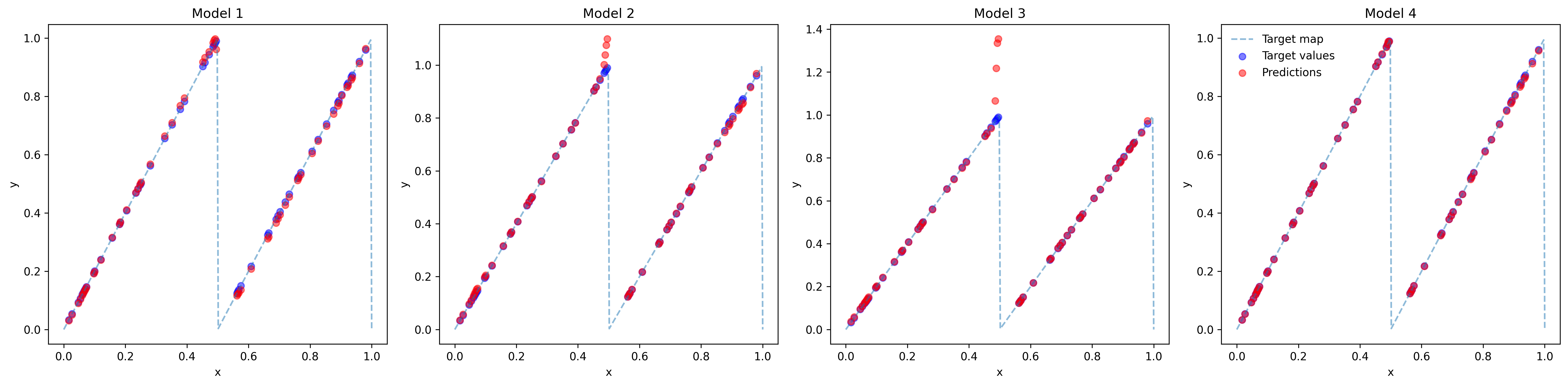}
        \caption{Doubling map}
        \label{fig:dya_pred}
    \end{subfigure}

    \begin{subfigure}[b]{\textwidth}
        \includegraphics[width=\textwidth]{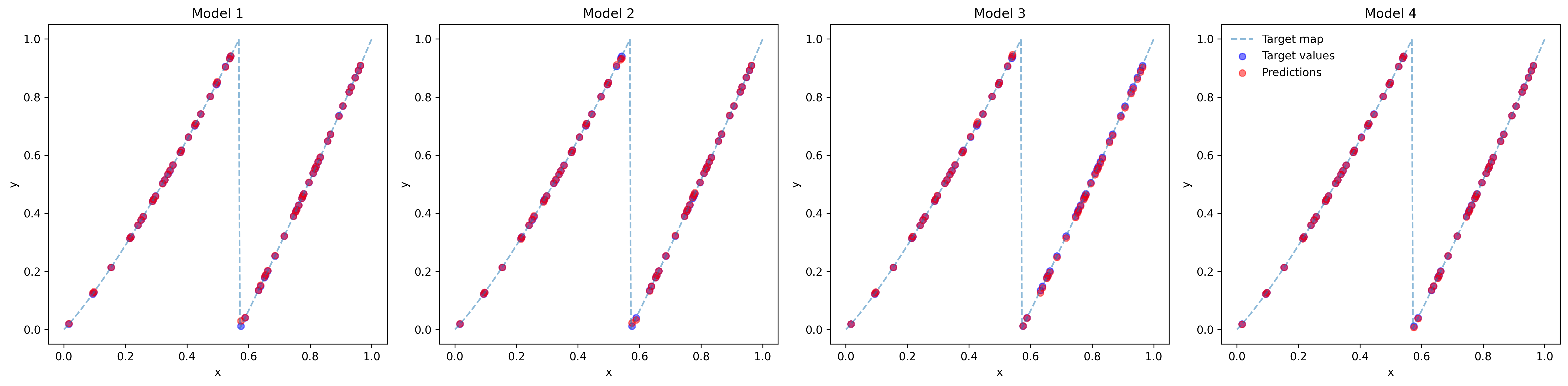}
        \caption{Pomeau-Manneville map}
        \label{fig:pm_pred_map}
    \end{subfigure}
    
    \caption{The prediction of the test dataset using deep learning models in \ref{auto}, including AENN (Models 1 and 2), FNN (Model 3), and PINN (Model 4), is based on Logistic map (A), Custom map (B), Katsura-Fukuda map (C), Doubling map (D) and Pomeau-Manneville map (E). The models used sequences of 300 in length, with 80\% of the data allocated for training and 20\% for testing.}
    
    \label{fig:4models}
\end{figure}

The prediction errors in Table \ref{tab:loss} indicate that the proposed model (Model 1) provides the lowest prediction error for continuous maps including logistic map, custom map, and Katsura-Fukuda map, while the PINN (Model 4) provide the lowest error for piecewise maps including doubling map and Pomeau-Manneville map. The comparison of predictions from different models against the actual values is shown in Figure \ref{fig:4models}. 

The proposed model successfully captures the dynamics of continuous maps compared to the FNN (Model 3) and PINN (Model 4). Our model demonstrates improved accuracy in capturing the dynamics of continuous maps compared to piecewise maps. The predictions for piecewise maps are less accurate than that for continuous maps, as illustrated in Figures \ref{fig:4models}. The model performs well in smooth regions but struggles in the discontinuous parts of the maps, where the functions change sharply. In addition, the baseline autoencoder (Model 2), which incorporates the logistic map in its latent space, and the standard neural network FNN exhibit larger deviations in these discontinuous regions. Notably, the errors are asymmetric, meaning the models do not mispredict in a consistent manner—some values are overestimated while others are underestimated. The mathematical structure of this discontinuous and intermittent behaviour, which is related to transitions between regular and chaotic phases, complicates the analysis of these maps, making them harder to study with networks that are more effective with smooth maps. We observe low prediction errors in piecewise maps for the PINN, which incorporates the informed dynamics map equation into its loss function to reduce prediction errors, but it does not yield the lowest prediction error for continuous maps. This indicates that although the PINN was informed by the precise physics, it does not guarantee accurate predictions. It suggests that incorporating meaningful features related to dynamical systems significantly enhances the performance of the neural network in specific ways.

\begin{figure}[htbp]
    \centering
    \begin{subfigure}[b]{\textwidth}
        \includegraphics[width=\textwidth]{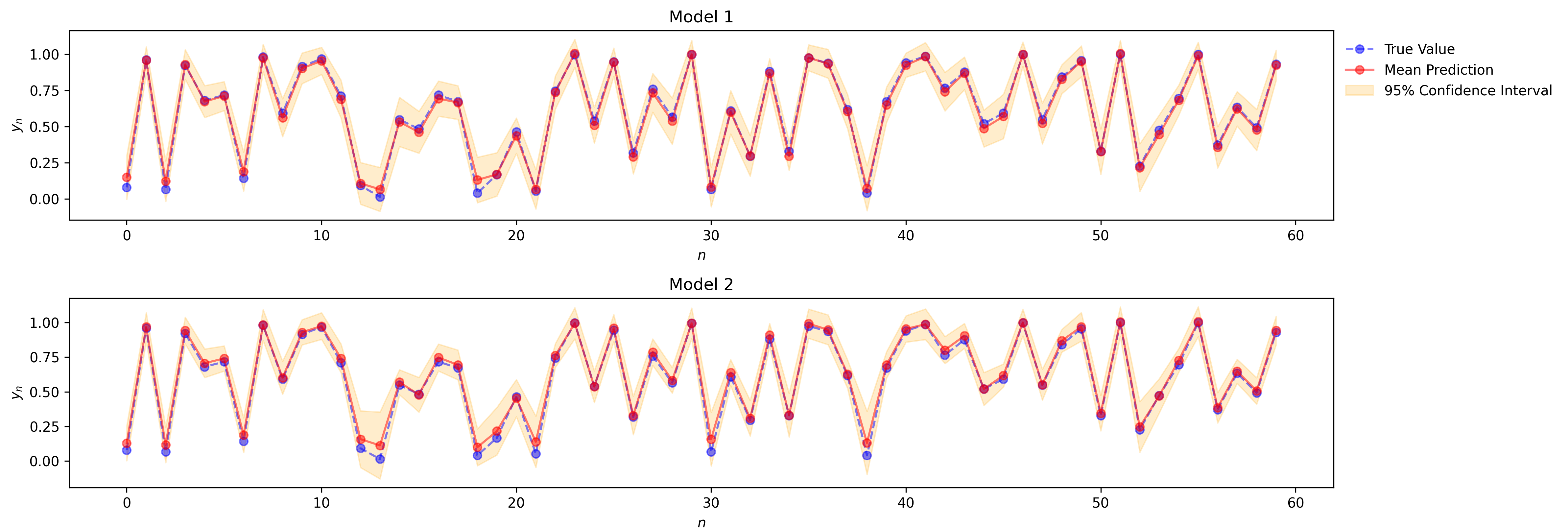}
        \caption{MC Dropout method}
        \label{fig:logistic_pred_model1}
    \end{subfigure}
    \hfill
    \begin{subfigure}[b]{\textwidth}
        \includegraphics[width=\textwidth]{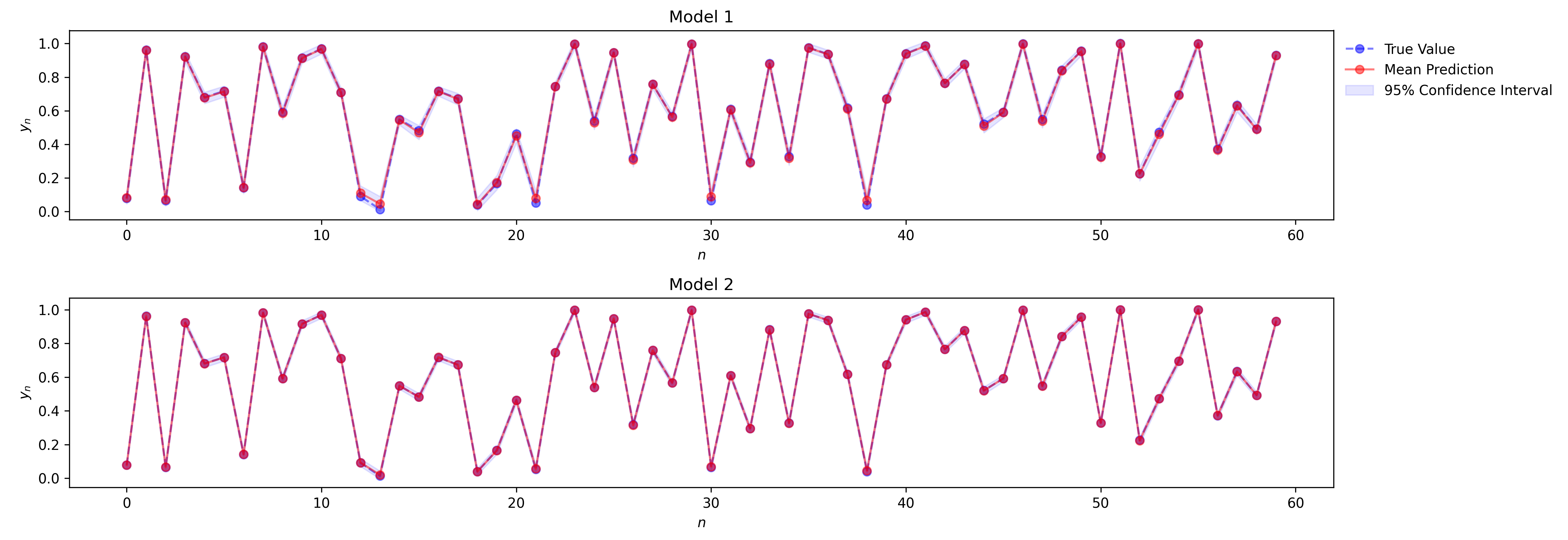}
        \caption{Ensemble method}
        \label{fig:logistic_pred_model2}
    \end{subfigure}
    
    \caption{Logistic map prediction results and 95\% confidence interval computed from MC dropout (top) and Ensemble methods (bottom) of AENN embedding conjugacy between tent map and logistic map in the latent space (Model 1) and AENN embedding the logistic map in the latent space (Model 2)}
    \label{fig:logistic_ci_prediction}
\end{figure}

\begin{figure}[htbp]

    \centering
    
    \begin{subfigure}[b]{\textwidth}
        \includegraphics[width=\textwidth]{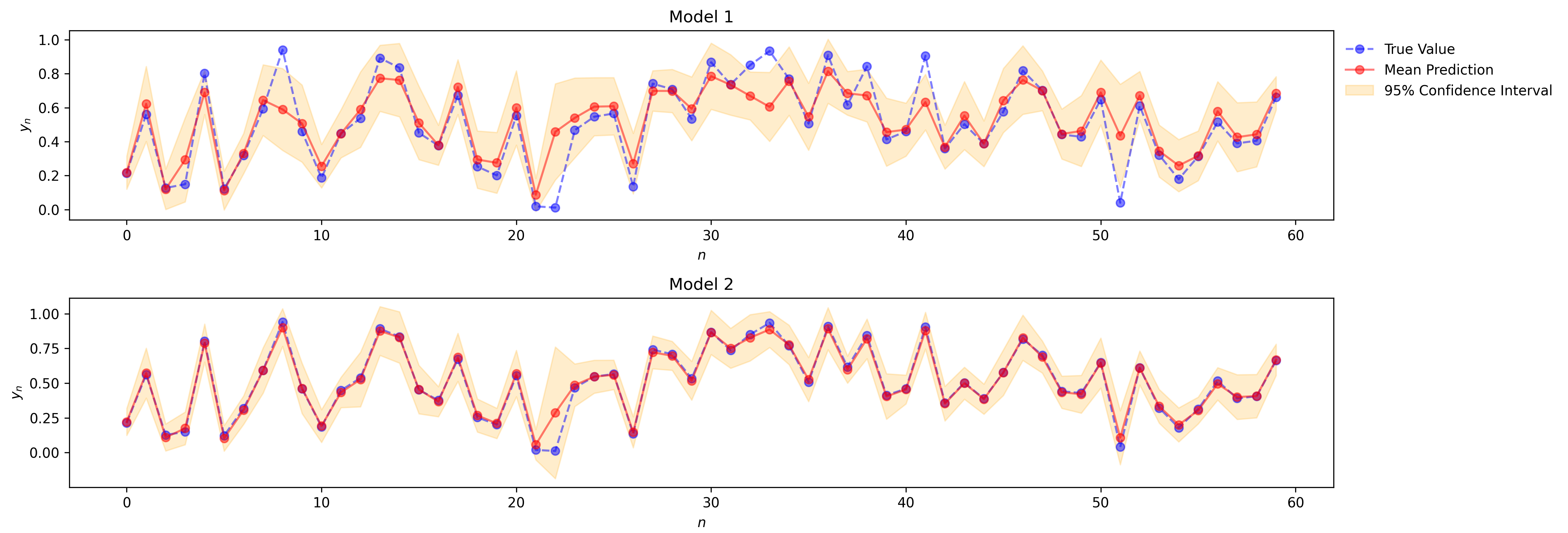}
        \caption{MC Dropout method}
        \label{fig:pom_pred_model1}
    \end{subfigure}
    \vspace{0.5cm}
    \begin{subfigure}[b]{\textwidth}
        \includegraphics[width=\textwidth]{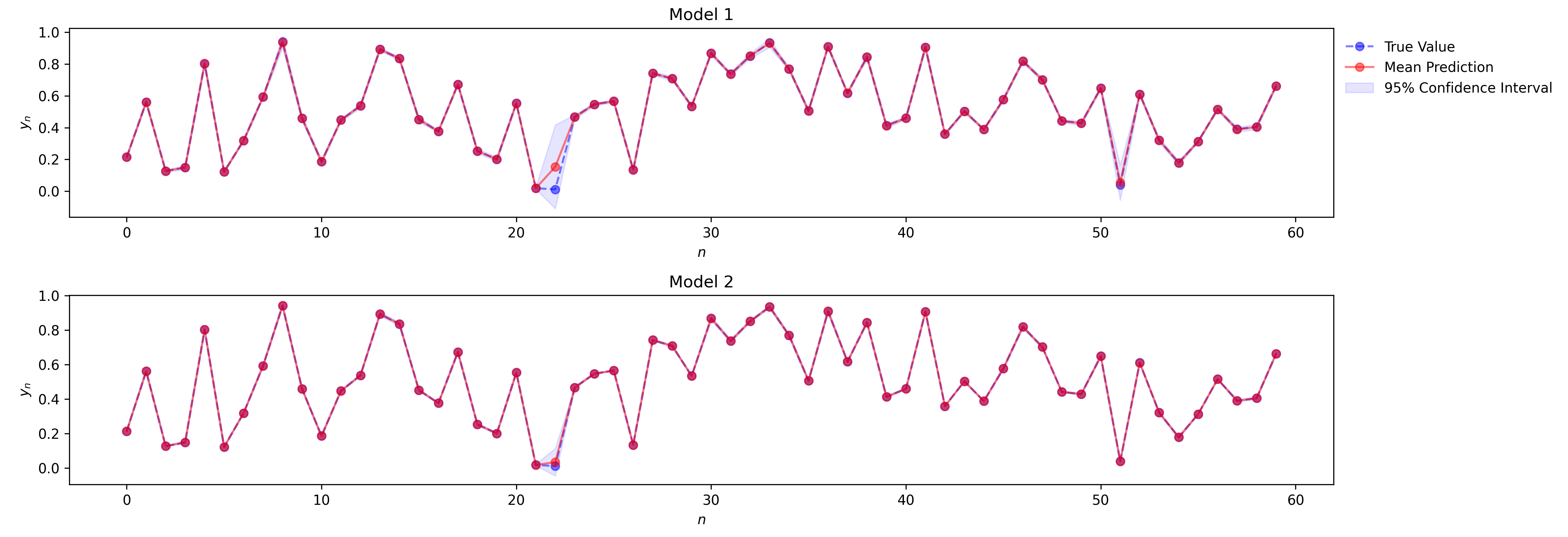}
        \caption{Ensemble method}
        \label{fig:pom_pred_model2}
    \end{subfigure}
    
    \caption{Pommeau-Manneville map prediction results and 95\% confidence interval computed from MC dropout (top) and Ensemble methods (bottom) of AENN embedding conjugacy between tent map and logistic map in the latent space (Model 1) and AENN embedding the logistic map in the latent space (Model 2)}
    \label{fig:pom_ci_prediction}
\end{figure}

In addition, we evaluate the uncertainties produced by our model using a 95\% confidence interval, which includes the mean and standard deviation. We compare the results to those from the baseline autoencoder (Model 2). The results are displayed in Figures \ref{fig:logistic_ci_prediction} and \ref{fig:pom_ci_prediction}. The MC dropout method, with a dropout probability of 0.2, indicates larger uncertainties across all maps. This observation suggests that the hidden dimension in the hidden layers significantly impacts both prediction accuracy and uncertainty estimation. Hence, balancing the hidden dimension is important when using the models. The confidence interval from the ensemble method, which incorporates a varying choice of initial weights for the networks, shows that the initialisation of weights produces relatively small uncertainties in the model's predictions.

\begin{figure}[htbp]
    \centering
    \begin{subfigure}[b]{\textwidth}
        \includegraphics[width=\textwidth]{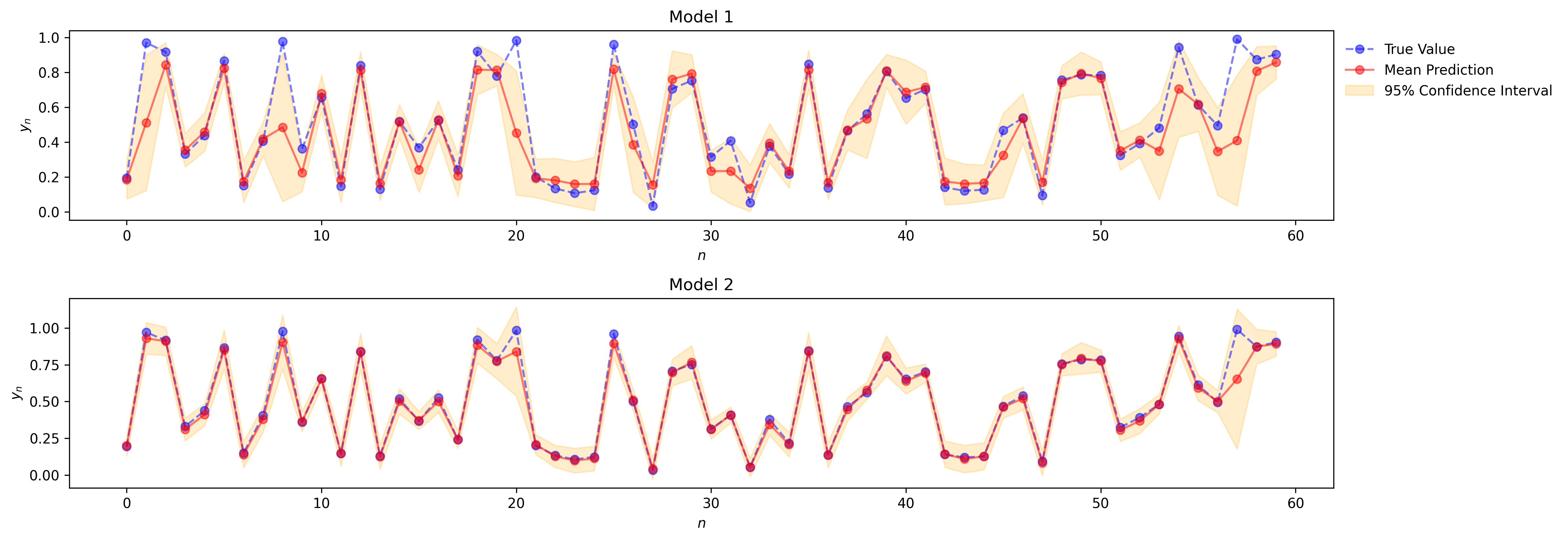}
        \caption{MC Dropout method}
        \label{fig:dy_pred_model1}
    \end{subfigure}
    \hfill
    \begin{subfigure}[b]{\textwidth}
        \includegraphics[width=\textwidth]{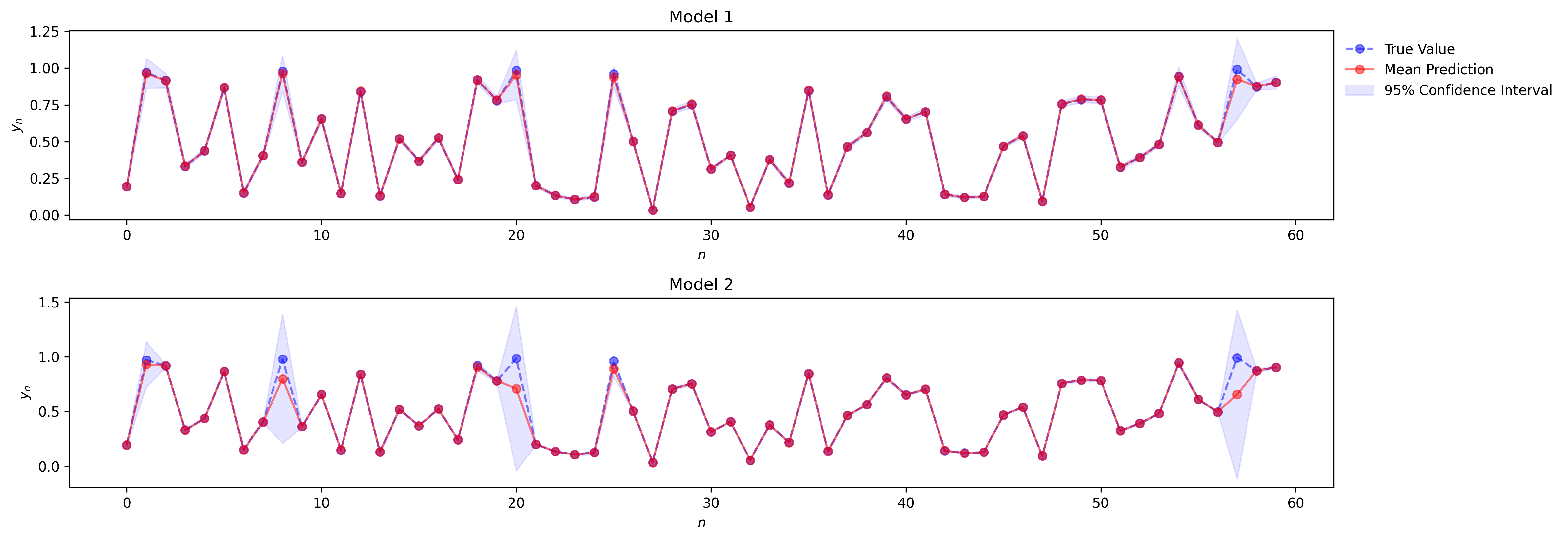}
        \caption{Ensemble method}
        \label{fig:dy_pred_model2}
    \end{subfigure}
    
    \caption{Doubling map prediction results and 95\% confidence interval computed from MC dropout (A) and Ensemble methods (B) of AENN embedding conjugacy between tent map and logistic map in the latent space (Model 1) and AENN embedding the logistic map in the latent space (Model 2)}
    \label{fig:dy_ci_prediction}
\end{figure}

Considering the conjugacy between target maps and the logistic map, the continuous maps with known conjugacy to the logistic map reveal that the proposed model performs as well as, or better than, the baseline autoencoder (Model 2). However, the performance of the baseline model is dependent on the choice of initial learnable parameters, which can lead to instability across different selections. In cases where there is no conjugacy or only semi-conjugacy exists, the autoencoder models that incorporate the conjugacy of the logistic map tend to produce higher errors, which error accumulation occurs at points that induce semiconjugacies, such as points of discontinuity. This was particularly evident with the doubling map, which is semiconjugate to the logistic map, and the Pomeau-Manneville map, which is neither conjugate nor semiconjugate to it. Our findings indicate less accurate predictions for the piecewise maps, as illustrated in Figure \ref{fig:4models}, and larger confidence intervals in Figures \ref{fig:dy_ci_prediction} and \ref{fig:pom_ci_prediction}. 



Overall, the findings suggest that the proposed model effectively analyses the behaviour of dynamical systems. The coordinate transformation using conjugacy between the chaotic tent map and the logistic map, enables effective training and approximation of the autoencoder compared to directly embedding the logistic map in the network. Although the model encounters challenges in accurately capturing discontinuities and sharp transitions in piecewise maps, it demonstrates significant improvements in deep learning performance.

\begin{landscape}
    \begin{table}[htbp]
    \centering
    \scriptsize
    \begin{tabular}{lccccccccc}
            \hline
            \multirow{3}{*}{Map} & \multirow{3}{*}{Parameter} & \multirow{3}{*}{Model 1}	&	\multicolumn{5}{c}{Model 2}	&	 \multirow{3}{*}{Model 3}	& \multirow{3}{*}{Model	4}	\\
	       &		&		& $c_1=3.0$	& $c_1=3.1$	& $c_1=3.5$ & $c_1=3.9$ & $c_1=4.0$ &&\\
    	   &		&		& $c_2=-3.0$ & $c_2=-3.1$ & $c_2=-3.5$ &	$c_2=-3.9$ & $c_2=-4.0$	&&\\
            \hline
            Logistic &	4.00	&	\cellcolor{blue!20} 1.528E-06	&	-	&	4.959E-06	&	9.459E-06	&	4.180E-06	&	-	&	2.106E-05	&	1.004E-05	\\
	               &	3.90	&	\cellcolor{blue!20} 1.472E-06	&	-	&	4.645E-06	&	2.462E-06	&	4.177E-06	&	-	&	5.641E-04	&	1.066E-05	\\
	               &	3.57	&	\cellcolor{blue!20} 1.045E-06	&	-	&	6.869E-06	&	7.005E-06	&	1.046E-05	&	-	&	6.888E-05	&	1.614E-03	\\
            
            Custom	&	-	&	\cellcolor{blue!20} 4.454E-06	&	-	&	5.329E-05	&	9.955E-06	&	8.267E-05	&	-	&	4.361E-05	&	8.130E-05	\\
            
            Katsura-Fukuta	&	0.50	&	\cellcolor{blue!20} 1.367E-06	&	-	&	4.959E-06	&	9.459E-06	&	3.133E-05	&	-	&	3.806E-04	&	8.377E-05	\\
            \hline 
            Doubling	&	-	&	9.089E-05	&	-	&	1.466E-03	&	4.742E-03	&	4.519E-04	&	-	&	5.415E-03	&	\cellcolor{blue!20} 7.066E-06 \\

            Pomeau-Manneville	&	1.50	&	7.444E-06	&	-	&	7.243E-06	&	7.395E-06	&	7.273E-06	&	-	&	1.381E-05	&	\cellcolor{blue!20} 1.222E-06	\\
            \hline
    \end{tabular}
    \caption{The prediction errors from different networks, '-' indicates a vanishing gradient. Hyperparameters were fixed across all models for computational comparison are shown in Table  \ref{tab:hyperparameter_search}.}\label{tab:loss}
\end{table}
\end{landscape}

\subsection{Limitations}
Here we have found that the piecewise linearity of the tent map enhances our network's ability to capture discontinuities when analysing piecewise maps, whilst the conjugacy function enables the network to learn continuous maps more efficiently. To achieve this, the latent space in our proposed network is strictly governed by the conjugacy between the chaotic tent map and the logistic map, and some assumed conjugacy between the logistic map and our target map of interest. This assumption makes the network less flexible for specific target dynamical maps that would otherwise violate this conjugate assumption. This problem is illustrated experimentally with less accurate predictions of piecewise maps than anticipated. Further improvements of this model may be achieved by reconstructing the latent space through various conjugacy functions derived from alternative maps, potentially resulting in more precise forecasting of the dynamics associated with piecewise maps. Nevertheless, it is important to note that finding a topological conjugacy function between two dynamical systems is usually challenging, and a general method is not available for arbitrary systems. Several studies highlight the challenges of conjugacy. A homeomorphism is required for topological conjugacy; however, its existence is not guaranteed. Determining a homeomorphism function between two maps remains difficult due to the specific structures, dimensions, and topological properties of the underlying spaces of dynamical systems \cite{bollt2010comparing}. The works of \cite{bollt2010comparing, skufca2008concept, skufca2007relaxing} discuss the concept of conjugacy functions and present a relaxation of this function based on the ideas of two maps being \textit{mostly conjugate} and their relationships via defective homeomorphisms. Further developments in this space, in particular incorporating results from the work of \cite{bollt2010comparing, skufca2008concept, skufca2007relaxing} noted above, could be made to improve accuracy or provide lower bounds on the error achievable by conjugacy autoencoder networks to learn chaotic maps.

\noindent\textbf{Data Availability.} All data and code used in this manuscript can be downloaded under the open access Creative Commons license from https://doi.org/10.48610/cc06c44 .

\noindent\textbf{Conflict of Interest Disclosure.} Meagan Carney reports travel was provided by University of Vienna Research Platform Erwin Schr\"{o}dinger International Institute for Mathematics and Physics. Hongkun Zhang reports travel was provided by The University of Queensland School of Mathematics and Physics. Ruethaichanok Kardkasem reports travel was provided by Great Bay University. If there are other authors, they declare that they have no known competing financial interests or personal relationships that could have appeared to influence the work reported in this paper.

\clearpage
\section{Appendix}

\begin{table}[htbp]
    
    \centering

    \begin{adjustbox}{width=\textwidth}
    
    \begin{tabular}{lrrrrrrr}
        \hline
         Map                & Parameter & Batch size & Layer Width & Layers In & Layers Out & Epoch & Learning rate\\
         \hline
         Logistic           & 4     &   64  & 256 & 1 & 1 & 1000& 0.005 \\
         Custom             & -     &   64  & 256 & 1 & 1 & 1000& 0.005 \\
         Katsura-Fukuda     & 0.5   &   64  & 256 & 1 & 1 & 1000& 0.005 \\
         Doubling             & -   &   64  & 256 & 3 & 3 & 2000 & 0.001 \\
         Pomeau-Manneville  & 1.5   &   64 & 64 & 3 & 3 & 2000& 0.001 \\
         \hline
    \end{tabular}

    \end{adjustbox}
    \caption{The hyperparameters used for training the neural networks are shown above. Note that the term "layer in" refers to the number of network layers used to build $h$, while "layer out" refers to the number of layers used to construct $h^{-1}$.  Although these two values are set to be the same, this is not a strict requirement.}
    \label{tab:hyperparameter_search}
\end{table}

\begin{figure}[htbp]
    \centering
    \begin{subfigure}[b]{\textwidth}
        \includegraphics[width=\textwidth]{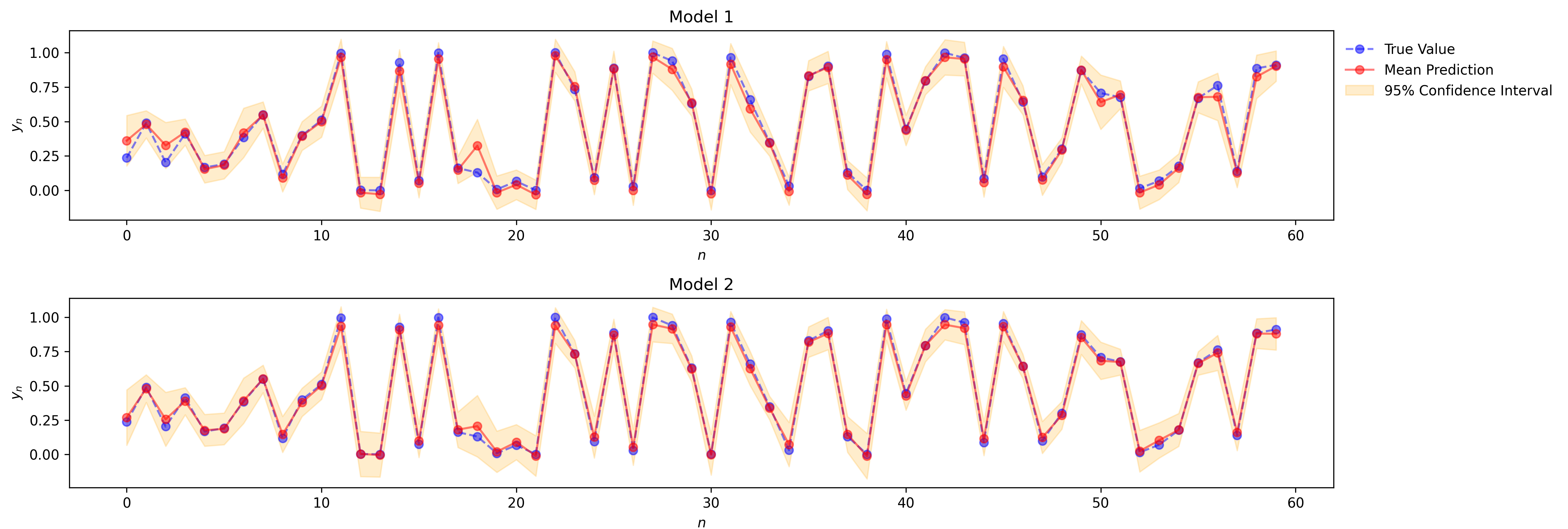}
        \caption{MC Dropout method}
        \label{fig:custom_pred_model1}
    \end{subfigure}
    \hfill
    \begin{subfigure}[b]{\textwidth}
        \includegraphics[width=\textwidth]{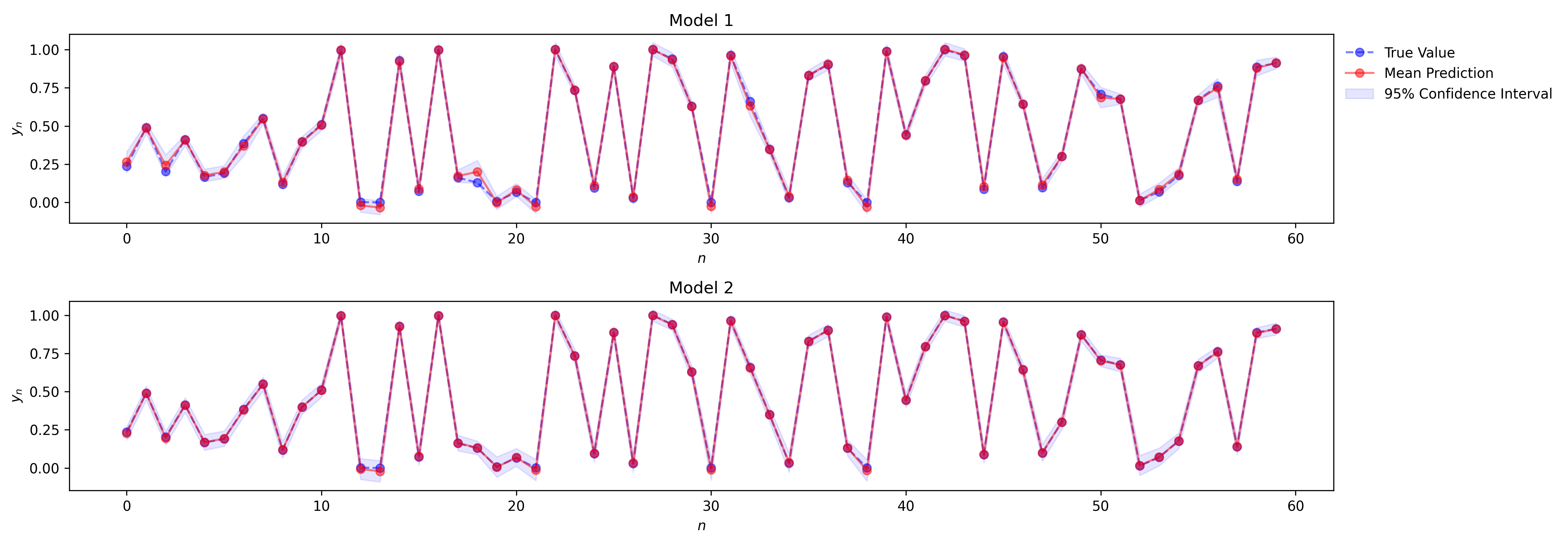}
        \caption{Ensemble method}
        \label{fig:custom_pred_model2}
    \end{subfigure}
    
    \caption{Custom map prediction results and 95\% confidence interval computed from MC dropout (top) and Ensemble methods (bottom) of AENN embedding conjugacy between tent map and logistic map in the latent space (Model 1) and AENN embedding the logistic map in the latent space (Model 2)}
    \label{fig:custom_ci_prediction}
\end{figure}

\begin{figure}[htbp]
    \centering
    \begin{subfigure}[b]{\textwidth}
        \includegraphics[width=\textwidth]{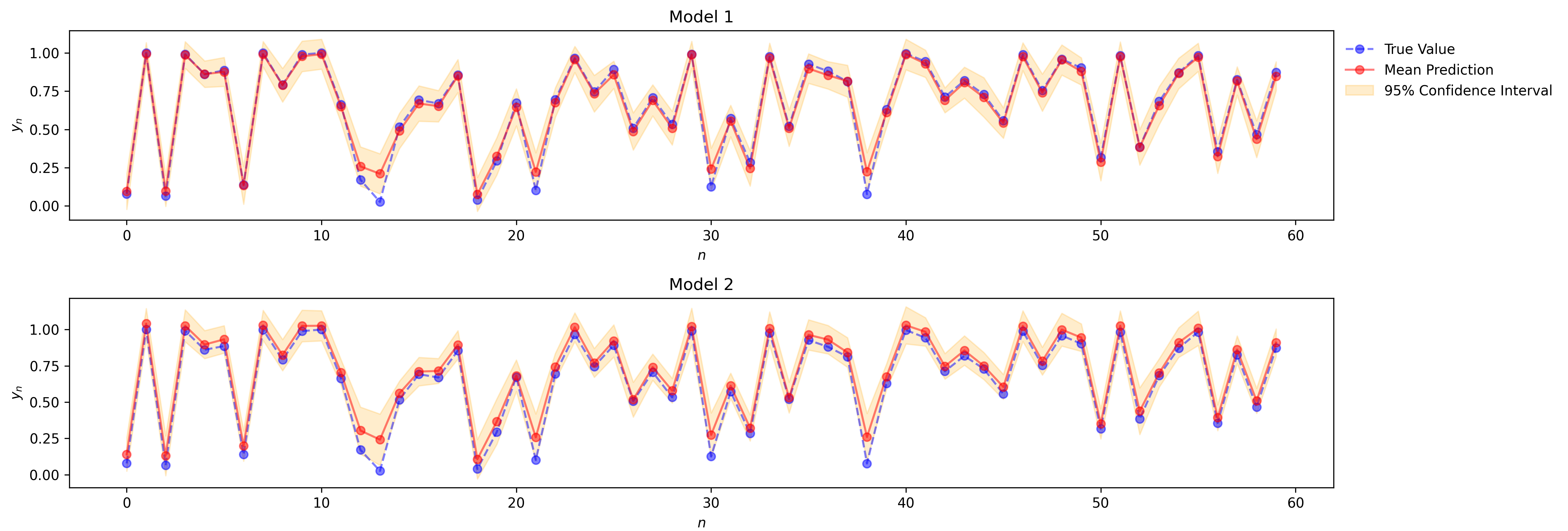}
        \caption{MC Dropout method}
        \label{fig:ks_pred_model1}
    \end{subfigure}
    \hfill
    \begin{subfigure}[b]{\textwidth}
        \includegraphics[width=\textwidth]{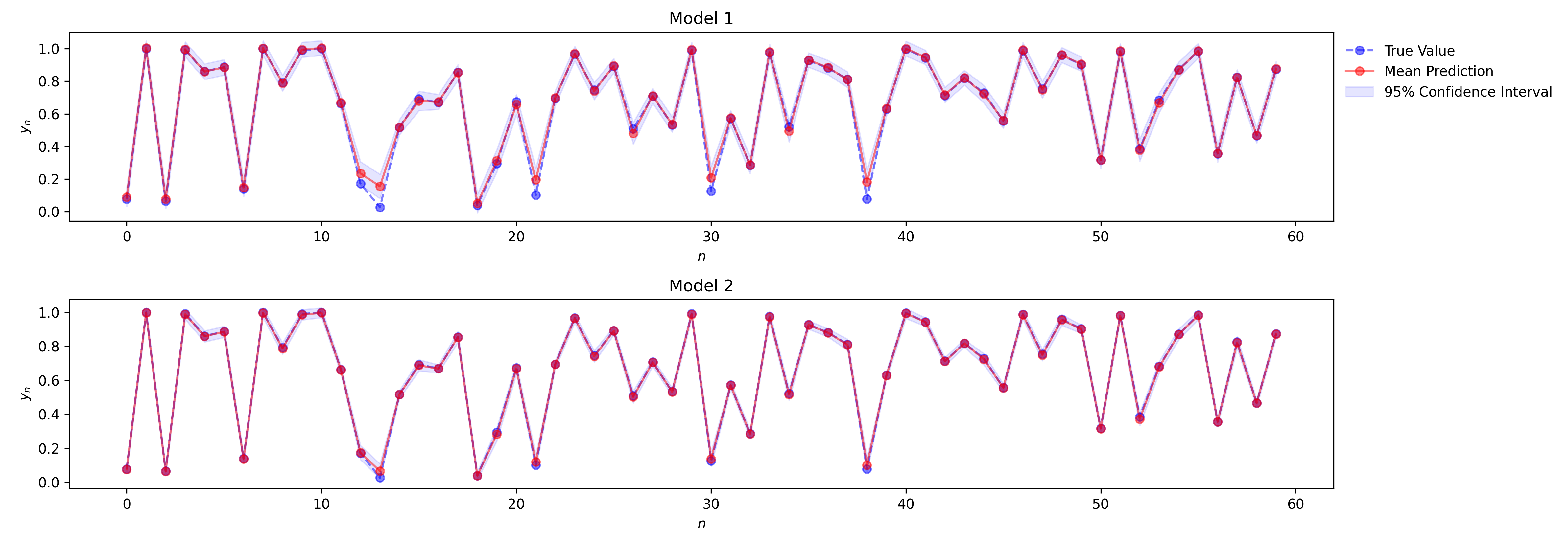}
        \caption{Ensemble method}
        \label{fig:ks_pred_model2}
    \end{subfigure}
    
    \caption{Katsura-Fukuda map prediction results and 95\% confidence interval computed from MC dropout (A) and Ensemble methods (B) of AENN embedding conjugacy  between tent map and logistic map in the latent space (Model 1) and AENN embedding the logistic map in the latent space (Model 2)}
    \label{fig:ks_ci_prediction}
\end{figure}

\clearpage
\bibliographystyle{unsrtnat}
\bibliography{reference}

\end{document}